# OPTIMAL AND FAST DETECTION OF SPATIAL CLUSTERS WITH SCAN STATISTICS[1]

By Guenther Walther

*Stanford University*

We consider the detection of multivariate spatial clusters in the Bernoulli model with $N$ locations, where the design distribution has weakly dependent marginals. The locations are scanned with a rectangular window with sides parallel to the axes and with varying sizes and aspect ratios. Multivariate scan statistics pose a statistical problem due to the multiple testing over many scan windows, as well as a computational problem because statistics have to be evaluated on many windows. This paper introduces methodology that leads to both statistically optimal inference and computationally efficient algorithms. The main difference to the traditional calibration of scan statistics is the concept of grouping scan windows according to their sizes, and then applying different critical values to different groups. It is shown that this calibration of the scan statistic results in optimal inference for spatial clusters on both small scales and on large scales, as well as in the case where the cluster lives on one of the marginals. Methodology is introduced that allows for an efficient approximation of the set of all rectangles while still guaranteeing the statistical optimality results described above. It is shown that the resulting scan statistic has a computational complexity that is almost linear in $N$.

**1. Introduction and overview of results.** Spatial scan statistics are used to detect clusters in spatial data and are widely used, for example, in epidemiology, biosurveillance and astronomy. In this paper, we consider the Bernoulli model in $\mathbf{R}^2$, which is important in many of the above applications, see, for example, Kulldorff (1999). All of the results in this paper can easily be extended to higher dimensions, but focusing on the important two-dimensional case simplifies the exposition and the notation. The Bernoulli model states that there are $N$ locations in $\mathbf{R}^2$, and that each location has a

Received October 2008; revised June 2009.
[1]Supported by NSF Grant DMS-05-05682 and NIH Grant 1R21AI069980.

label associated with it that takes on one of two possible outcomes, say 0 and 1. Conditional on the locations, the values of these labels are realizations of independent Bernoulli random variables with parameter $p$ at all locations in a certain set $R$ and parameter $q$ at all locations in $R^c$.

The null hypothesis is

$$H_0 : p = q$$

and the alternative hypothesis is $H_1 : q < p$ for some unknown set $R$ out of a certain class of sets, which in this paper is taken to be the class of rectangles with sides parallel to the axes and with arbitrary sizes and aspect ratios.

EXAMPLE 1. Each location represents the spatial location of a person which is either healthy (label 0) or diseased (label 1). If $q < p$, then $R$ represents the local area of a disease outbreak. The task is to detect such disease clusters $R$ where the disease density is significantly higher than the population density; see, for example, Kulldorff (1999) and the references given there.

EXAMPLE 2. A flow cytometer measures various numerical characteristics of each of a large number of cells. Thus, each cell can be identified with a point in Euclidean space. One task in the analysis flow cytometry data is to describe local regions where two cell distributions are different. Roederer and Hardy (2001) and Roederer et al. (2001) model such regions as rectangles with sides parallel to the axes. The reason for this is that such sets are easy to interpret and easy to implement on instruments for further processing. Also, sometimes the difference effect lives on a lower-dimensional subspace, which is a special case of axis-parallel rectangles. Suppose one has a sample drawn i.i.d. from a distribution $G$, and a second sample drawn i.i.d. from a distribution $H$. Label each observation in the first sample with a 0 and each observation in the second sample with a 1. If $G = H$, then the problem is left invariant under permutations of the labels. It will be seen shortly that thus the inference for Example 2 can proceed identically to that for Example 1.

Further examples are given, for example, in Kulldorff (1999). There is a large body of work on univariate scan statistics; see, for example, the references in Glaz and Balakrishnan (1999) and in Glaz, Naus and Wallenstein (2001), but the multivariate case is much less well developed. One reason is that computational issues play a prominent role when multivariate scan windows need to be evaluated at possibly many locations. Some references for the multivariate case that are relevant for the problem studied here are Naus (1965), Loader (1991), Chen and Glaz (1996), Alm (1997), Anderson and Titterington (1997), Kulldorff (1997, 1999), Naiman and Priebe (2001), and on the computational aspect, Neill and Moore (2004a, 2004b). A recent reference for the multivariate two-sample problem is Rohde (2009), who constructs regions of significant difference based on nearest-neighbor statistics.

To simplify the exposition, we will assume that the distribution $F$ of the locations is continuous and has independent marginals. All of the main results continue to hold if the marginals are weakly dependent, for example, if they are $\psi$-mixing. But a completely general design of the $N$ locations requires some modifications, which will be reported elsewhere. It is also convenient to only consider rectangles $R$ with $F(R) \leq 1/8$, which is an innocuous restriction for most problems.

In our analysis, we will condition on the sample size $N$ and on the $N$ locations. Under the null hypothesis, the problem is then left invariant under permutations of the labels, and exact finite sample significance statements for rectangles $R$ can be obtained by a permutation test. There are two major problems associated with such an inference: as the class of rectangles is large, a statistical problem arises in the form of multiple testing, and a computational problem arises due to the need to evaluate test statistics on many rectangles. This paper introduces methodology that leads to both statistically optimal inference and computationally efficient algorithms.

The conventional definition of a scan statistic is

$$\max_{R \in \mathcal{R}} T(R), \tag{1}$$

where $\mathcal{R}$ is a set of scan windows such as the set of rectangles described above, and $T$ is a standardized test statistic that is evaluated locally for each scan window. Critical values are then derived for this overall maximum. In this paper, we propose to use size-dependent critical values obtained by grouping windows according to their size: all windows that contain between $2^{-\ell-1}N$ and $2^{-\ell}N$ locations are grouped into one block, $\ell \geq 3$. Then we use different critical values for different blocks as proposed by Rufibach and Walther (2009) in a certain univariate context. The heuristic motivation for this approach is the following: there are of the order $N$ disjoint windows containing a small number of locations each. As the corresponding local statistics $T$ will be roughly independent, the maximum over the small windows will behave like the maximum of $N$ i.i.d. random variables. This will tend to be stochastically much larger than the maximum over windows of size $N/8$ (say), which will roughly behave like the maximum of 8 of these i.i.d. random variables. Thus, the distribution of the conventional scan statistic (1) will be dominated by the small windows, with a corresponding loss of power for larger windows. Grouping windows according to their size and employing size-dependent critical values allows to remedy this effect. Indeed, it will be shown below that this methodology allows for the following large sample results:

- If the effect $q < p$ lives on a small rectangle, then the blocked scan statistic is essentially as powerful as any test can possibly be.
- If the effect $q < p$ lives on a large rectangle, then the blocked scan statistic is again optimal, even in comparison to tests that are allowed to use a priori knowledge of the correct window size. That is, scanning with different window sizes does not result in a significant penalty.
- If the effect $q < p$ lives on one of the two marginals, then the above optimality results still hold in the one-dimensional framework. That is, scanning with two-dimensional rectangles does not result in a significant penalty even if it is known a priori that the effect lives on a univariate marginal.

We will give heuristic explanations of these results as well as rigorous mathematical statements. These results use a concentration inequality for the hypergeometric distribution which may be of independent interest. The optimality results require the use of size-dependent critical values, and it appears that such methodology has not been used before for scan statistics. In the setting of univariate function estimation in the Gaussian White Noise model, Dümbgen and Spokoiny (2001) employed a scale-dependent penalty term. It is not clear how a useful penalty term can be derived for the problem under consideration here. Also, the univariate results in Rufibach and Walther (2009) suggest that the block procedure yields a better finite sample performance for relevant sample sizes.

The construction of efficient algorithms, and also the particular proof of the above optimality results, requires an economical approximation of the set of all rectangles. We prove an approximation theorem that allows for an adequate approximation of the set of all rectangles by $O(N \log^2 N)$ rectangles. By comparison, there are $O(N^4)$ rectangles that contain different subsets of the $N$ locations. As a consequence, it will be shown that the blocked scan statistic can be implemented with a computational complexity that is almost linear in $N$.

It will be seen that there is a close connection between the computational approximation scheme and the statistical inference, with the grouping of rectangles according to their size being a central theme in each case.

**2. The blocked scan statistic.** Kulldorff (1997) derives the log-likelihood ratio statistic for a given scan window $R$ as

$$T(R) = n\left(\hat{p}\log\frac{\hat{p}}{\bar{p}} + (1-\hat{p})\log\frac{1-\hat{p}}{1-\bar{p}}\right)$$
$$+ (N-n)\left(\hat{q}\log\frac{\hat{q}}{\bar{p}} + (1-\hat{q})\log\frac{1-\hat{q}}{1-\bar{p}}\right),$$

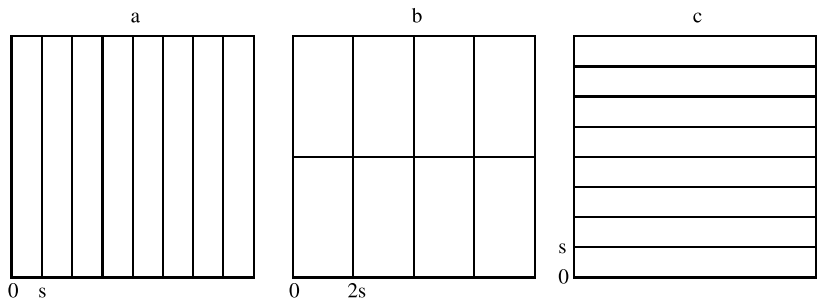

Fig. 1. *Constructing an approximating set of rectangles. Units on the axes are with respect to the marginal distributions $F_X$ and $F_Y$.*

if $\hat{q} \leq \hat{p}$, and $T(R) = 0$, otherwise. Here, $n := \#R$ is the number of locations in $R$, $\overline{p}$ is the overall proportion of 1's, and $\hat{p}$ and $\hat{q}$ are the proportion of 1's in $R$ and in $R^c$, respectively. Despite its cumbersome form, this statistic has been widely adapted in the computer science literature; see, for example, Neill and Moore (2004a, 2004b). The concentration inequality given in Theorem 4 shows that this transformation of $\hat{q}$ and $\hat{p}$ has the benefit of a clean tail behavior.

$T(R)$ is zero if $\hat{p} = \hat{q}$ and positive if $\hat{q} < \hat{p}$. We restrict ourselves to the alternative hypothesis $q < p$ for notational simplicity and because this case is most relevant for applications, but all the following results continue to hold for the alternative $p \neq q$ after a simple modification of the definition of $T(R)$.

We will evaluate $T$ on a set of rectangles that is a good approximation to the set $\mathcal{R} := \{\text{axis-parallel rectangles in } \mathbf{R}^2\}$. The following theorem has a constructive proof that shows how to construct an economical set of rectangles that approximate all rectangles in $\mathcal{R}$ whose size (as measured in terms of $F$) is about $s$, namely $\mathcal{R}(s) := \{R \in \mathcal{R} : s/2 < F(R) \leq s\}$.

THEOREM 1. *For every $s, \varepsilon \in (0,1)$, there exists $\mathcal{R}_{\text{app}}(s,\varepsilon) \subset \mathcal{R}$ such that:*

1. *For every $R \in \mathcal{R}(s)$, there exists $R' \in \mathcal{R}_{\text{app}}(s,\varepsilon)$ with $F(R \triangle R') \leq \varepsilon F(R)$.*
2. *$\#\mathcal{R}_{\text{app}}(s,\varepsilon) \leq Cs^{-1}\varepsilon^{-4}\log(2/s)$ for a universal constant $C$.*

The idea for the approximation scheme is depicted in Figure 1 and explained in Section 3. To construct an approximation for all of $\mathcal{R}$ we proceed as follows. First, note that $\mathcal{R} = \bigcup_{\ell=0}^{\infty} \mathcal{R}(2^{-\ell})$. Second, the construction of $\mathcal{R}_{\text{app}}(s,\varepsilon)$ depends on $F$, which is typically unknown. To obtain an approximation that depends on the observations only, we replace $F$ by the empirical measure $F_N$ in the construction of $\mathcal{R}_{\text{app}}(s,\varepsilon)$ and call the resulting set

$\mathcal{R}_{\mathrm{app},N}(s,\varepsilon)$. Then we define our approximating set as

$$\text{(2)} \qquad \mathcal{R}_{\mathrm{app},N} := \bigcup_{\ell=3}^{\lfloor \log_2(N/(2\log N)) \rfloor} \mathcal{R}_{\mathrm{app},N}(2^{-\ell}, \ell^{-1/2}).$$

The particular choice $\varepsilon = \ell^{-1/2}$ yields the optimality results given in Theorem 2 below. Thus, the smaller the rectangle, the finer the approximation relative to the size of the rectangle. Section 4 gives an algorithm that constructs $\mathcal{R}_{\mathrm{app},N}$. The precise result about this approximation is as follows.

COROLLARY 1. *There exists $\mathcal{R}_{\mathrm{app},N} \subset \mathcal{R}$ depending only on $F_N$ such that:*

1. *For every $R \in \mathcal{R}$ with $F(R) \in [\frac{2\log N}{N}, \frac{1}{8}]$, there exists $R' \in \mathcal{R}_{\mathrm{app},N}$ with $R' \subset R$ and $F(R \triangle R') \leq \frac{9}{8} \frac{F(R)}{\sqrt{\lfloor \log_2(1/F(R)) \rfloor}}$ with probability converging to 1 uniformly in $R \in \mathcal{R}$ and $F$.*
2. $\#\mathcal{R}_{\mathrm{app},N} \leq C'N \log^2 N$ *for a universal constant $C'$.*

By comparison, a naive enumeration of all rectangles that contain different subsets of the $N$ locations results in $O(N^4)$ rectangles and therefore is generally computationally infeasible. The algorithm in Section 4 computes the scan statistic $T$ over $\mathcal{R}_{\mathrm{app},N}$ in $O(N \log^4 N)$ steps, i.e., with a computation time that is almost linear in $N$.

We will call $\mathcal{R}_{\mathrm{app},N}(2^{-\ell}, \ldots)$ the $\ell$th block of rectangles. The idea for the statistical methodology is closely connected to this approximation scheme: as all the rectangles in the $\ell$th block have about the same size, we will assign to those rectangles the same critical value. Following the criterion given in Rufibach and Walther (2009), we set these critical values such that the significance level of the $\ell$th block decreases as $\sim \ell^{-2}$.

In more detail, let $\alpha \in (0,1)$ and define $q_\ell(\alpha)$ to be the $(1-\alpha)$-quantile of $\max_{R \in \ell\text{th block}} T(R)$ when the labels are permuted randomly. For notational convenience, we suppress the dependence of $q_\ell(\alpha)$ on the sample size $N$ and on $\overline{p}$. Let $\tilde{\alpha}$ be the largest number such that

$$\text{(3)} \qquad \mathbb{P}\left( \bigcup_{\ell=3}^{\lfloor \log_2(N/(2\log N)) \rfloor} \left\{ \max_{R \in \ell\text{th block}} T(R) > q_\ell\left(\frac{\tilde{\alpha}}{\ell^2}\right) \right\} \right) \leq \alpha.$$

By construction, one can then claim with guaranteed simultaneous finite sample confidence $1 - \alpha$ that $H_0$ is violated on every rectangle $R$ on which $T(R) > q_\ell(\tilde{\alpha}/\ell^2)$, where $\ell$ is the block index of $R$. As explained in Rufibach and Walther (2009), it is advantageous in practice to replace $\ell^2$ by, for example, $(10 + \ell)^2$, and all of the following results also apply for such a modification.

The $q_\ell$ can be readily simulated with a simple extension of the usual Monte Carlo technique for a permutation test: for each block, one records in a list the maximum for each Monte Carlo permutation of the labels. Then one can use a bisection method on the lists of sorted maxima to find $\tilde{\alpha}$.

**3. Optimality.** In the following, we consider a growing sample size $N$. The result below allows for a quite general situation where the rectangle $R_N$ may vary with $N$, likewise the probabilities of success $p_N$ in $R_N$ and $q_N$ in $R_N^c$, and the design distribution $F^N$. For simplicity, we denote probabilities under this model by $\mathbb{P}_N$. The key quantity for detecting $q < p$ on some rectangle $R$ turns out to be

$$D(F(R), p, q) := F(R)(1 - F(R))\frac{(p-q)^2}{p(1-q)}.$$

$D(F(R), p, q)$ increases, and hence the detection of $R$ becomes easier, if for fixed $F(R)$ and $q$ the difference $p - q$ increases, as one would expect. If $F(R)$ and $p - q$ are fixed, then $D(F(R), p, q)$ increases as $(p + q)/2$ moves away from $1/2$, i.e., detection is easier if the background probability $q$ is closer to 0 or if $p$ is closer to 1. Theorem 2 below quantifies when detection is possible and shows that the blocked scan statistic is optimal for detecting both small rectangles, that is, when $F^N(R_N) \to 0$, and large rectangles, that is, when $\liminf F^N(R_N) > 0$. An appropriate way to formulate these optimality results is via the asymptotic minimax framework, see, e.g., the investigation of univariate shape properties on small scales in the Gaussian white noise model in Dümbgen and Spokoiny (2001) and the results for small and large scales in the context of a univariate density in Dümbgen and Walther (2008).

THEOREM 2. (a) *Let $\{(F^N, R_N, p_N, q_N)\}$ be an arbitrary sequence of parameters with*

$$D(F^N(R_N), p_N, q_N) \geq (2 + \varepsilon_N)\frac{\log(1/F^N(R_N))}{N}$$

$$\text{where } \varepsilon_N \sqrt{\log \frac{1}{F^N(R_N)}} \to \infty.$$

*Then*

$\mathbb{P}_N(\text{the blocked scan statistic finds a significant rectangle } R \subset R_N) \to 1.$

(b1) *Let $\phi_N$ be any sequence of tests with asymptotic level $\alpha \in (0, 1)$ under $H_0$. For any prescribed sequence of continuous distributions $\{F^N\}$, there* *exists a sequence of parameters* $\{(R_N, p_N, q_N)\}$ *such that*

$$D(F^N(R_N), p_N, q_N) \geq (2 - \varepsilon_N) \frac{\log(1/F^N(R_N))}{N}$$

*with* $\varepsilon_N \downarrow 0, \varepsilon_N \sqrt{\log \frac{1}{F^N(R_N)}} \to \infty$, *and* $\overline{\lim}_N \mathbb{P}_N(\phi_N \text{ rejects}) \leq \alpha$.

*This result continues to hold if one also prescribes the values* $\{F^N(R_N)\}$ *and* $\{\varepsilon_N\}$, *provided that* $(\log N)^2/N \leq F^N(R_N) \to 0$ *and* $\varepsilon_N \sqrt{\log \frac{1}{F^N(R_N)}} \to \infty$.

(b2) *Let* $\{F^N\}$ *be any sequence of continuous distributions and* $\{R_N\}$ *any sequence of rectangles with* $F^N(R_N)(1 - F^N(R_N)) > 0$, *let* $b_N \in [0, NF^N(R_N)(1 - F^N(R_N)))$, *and let* $\phi_N$ *be any test with asymptotic level* $\alpha \in (0,1)$ *under the null hypothesis that the probability of success on* $R_N$ *equals that on* $R_N^c$. *If*

$$\mathbb{P}_N(\phi_N \text{ rejects}) \to 1$$

*for every sequence of parameters* $\{(p_N, q_N)\}$ *that satisfies* $D(F^N(R_N), p_N, q_N) \geq \frac{b_N}{N}$, *then necessarily* $b_N \to \infty$.

Parts (a) and (b1) show that in the case of small rectangles there is a cutoff at $D = 2\log \frac{1}{F^N(R_N)}/N$: if $D \geq (2 + \varepsilon_N)\log \frac{1}{F^N(R_N)}/N$ with $\varepsilon_N \to 0$ sufficiently slowly, then the blocked scan statistic will detect the rectangle $R_N$ with asymptotic power one. One the other hand, if $D$ is of the size $(2 - \varepsilon_N)\log \frac{1}{F^N(R_N)}/N$, then no test can exist that detects the rectangle with nontrivial asymptotic power. These two statements leave essentially no room for any other test to beat the blocked scan statistic for detecting small rectangles.

In the case of large rectangles, part (b2) states that any test $\phi_N$ can have asymptotic power 1 only if $ND \to \infty$. But under the latter condition, the blocked scan statistic also has asymptotic power 1 [because $ND \to \infty$ arbitrarily slowly is sufficient for the claim in (a) if $\frac{1}{F^N(R_N)}$ stays bounded]. Note that (b2) even allows the competing test $\phi_N$ to possess prior knowledge of the rectangle $R$.

These results clarify the tradeoff when using a scan statistic with varying window size. On the one hand, one can evidently gain substantial power by matching the window size with the extent of the effect. On the other hand, varying the window size incurs a multiple testing penalty. The above results show that this multiple testing penalty becomes negligible for large samples, provided one employs an appropriate calibration of the various window sizes such as the blocked scan statistic introduced here. An illustration will be given in Section 4.

Furthermore, there is no substantial multiple testing penalty for searching over multivariate rectangles when the effect lives on one of the marginals.

THEOREM 3. *Suppose that the $\{R_N\}$ are in fact intervals on one of the two axes. Then the conclusions of Theorem 2 continue to hold, even if the tests $\phi_N$ in* (b1) *and* (b2) *are allowed to use the prior knowledge about which axis the $\{R_N\}$ live on.*

A heuristic explanation of this result is as follows: Figure 1(a) depicts $\frac{1}{s}$ disjoint rectangles with content $s$. The rectangles of Figure 1(b) are obtained by doubling the width of certain rectangles in Figure 1(a) and then dividing the rectangle into two with a horizontal split. After $\log \frac{1}{s}$ iterations one obtains the rectangles of Figure 1(c). The idea of Theorem 1 is that in the case of independent or weakly dependent marginals the totality of these rectangles (after a refinement allowing, e.g., certain translations) constitutes an economical set of rectangles that approximates well the set of all rectangles with content $s$. The difficulty of the multiple testing problem depends essentially on the cardinality of this approximating set, as local statistics that pertain to rectangles with large overlap will be highly correlated and thus will not affect the multiple testing problem much. But the construction depicted in Figure 1 results in $\frac{1}{s}\log\frac{1}{s}$ rectangles, which up to the log term is of the same order as the $\frac{1}{s}$ rectangles in the univariate case of Figure 1(a). Thus, one expects that the multiple testing problem in this multivariate situation will not be significantly more difficult than in the univariate case.

The proof of Theorem 2 makes use of the following concentration inequality for the hypergeometric distribution.

THEOREM 4. *Let $X$ denote the number of red items among $n$ items drawn without replacement out of $N$ items of which $R$ are red. Then*

$$\mathbb{P}(X \geq x) \leq C(L(x)+2)\exp(-L(x)) \qquad \text{for } x > m := nR/N,$$
$$\mathbb{P}(X \leq x) \leq C(L(x)+2)\exp(-L(x)) \qquad \text{for } x < m$$

*and*

$$\mathbb{P}(L(X) \geq x) \leq 2C(x+2)\exp(-x) \qquad \text{for } x > 0,$$

*where $L(x) := n(\hat{p}\log\frac{\hat{p}}{\bar{p}} + (1-\hat{p})\log\frac{1-\hat{p}}{1-\bar{p}}) + (N-n)(\hat{q}\log\frac{\hat{q}}{\bar{p}} + (1-\hat{q})\log\frac{1-\hat{q}}{1-\bar{p}})$, $\bar{p} := \frac{R}{N}$, $\hat{p} := \frac{x}{n}$, $\hat{q} := \frac{R-x}{N-n}$, and $C := 2\exp\{\frac{13}{12\bar{p}(1-\bar{p})}(\frac{1}{n}+\frac{1}{N-n})\}$.*

This inequality compares to the classical concentration bound obtained from the Chernoff–Hoeffding theorem as follows. Hoeffding [(1963), Theorem 1 and Section 6] gives $\mathbb{P}(X \geq x) \leq \exp(-n(\hat{p}\log\frac{\hat{p}}{\bar{p}} + (1-\hat{p})\log\frac{1-\hat{p}}{1-\bar{p}}))$. A Taylor series expansion shows that for $\hat{p}$ near $\bar{p}$ the exponent behaves like $-n\frac{(\hat{p}-\bar{p})^2}{2\bar{p}(1-\bar{p})}$, whereas $-L(x) \approx -n\frac{(\hat{p}-\bar{p})^2}{(1-n/N)2\bar{p}(1-\bar{p})}$. Thus, Theorem 4 accounts for the variance correction factor $\frac{N-n}{N-1}$ for sampling without replacement.

Derbeko, El-Yaniv and Meir (2004) give an improvement to the Chernoff–Hoeffding bound that is weaker than that of Theorem 4. Rohde (2009) gives a Bernstein-type inequality for the hypergeometric distribution.

*Scanning on a grid and comparison with the algorithm of Neill and Moore.* Neill and Moore (2004a, 2004b) give an algorithm that runs in $O(N \log^2 N)$ steps for data that are binned on a $\sqrt{N} \times \sqrt{N}$ grid. That algorithm produces a rectangle $R$ that attains $\max_R T(R)$ by partitioning the grid into overlapping regions, bounding $\max T$ over subregions, and pruning regions which cannot contain the maximum. Thus, both the algorithm of Neill and Moore and the algorithm introduced here run in almost linear time; see Proposition 1 below. Both algorithms achieve this by using an approximation. The algorithm of Neill and Moore approximates the data by binning them on a grid, and then finds the exact maximum over all rectangles on the grid. In contrast, the methodology introduced here evaluates rectangles on the exact data, but approximates the set of all rectangles. The results of this section show that this algorithm results in a solution that is statistically optimal. It is an open problem whether the algorithm of Neill and Moore results in a solution that is statistically optimal, and how the grid has to be constructed to achieve this.

Consider now the related problem where one observes a Bernoulli random variable on each grid point of a $\sqrt{N} \times \sqrt{N}$ grid. Then the design distribution $F$ has independent marginals, but it is not continuous any more. Still, the methodology introduced in this paper can be readily adapted to this set-up and shown to be statistically optimal. The conclusions of Theorems 2 and 3 continue to hold. However, the condition on the size of $F(R)$ in Theorem 2(b1) now has to be set differently for the marginal effects considered in Theorem 3. In the multivariate case, (b1) allows rectangles $R$ as small as $F(R) \geq \log^2 N/N$ [and (a) even allows detection if $F(R) \geq 2\log N/N$], which results in a detection threshold of about $2\log N/N$ for rectangles with these sizes. But any nonempty marginal interval $R$ necessarily satisfies $F(R) \geq \sqrt{N}/N$ due to the nature of the grid, which results in a detection threshold of about $\log N/N$ for the smallest detectable marginal intervals.

*Controlling* $\max_R T(R)$. The cardinality of the approximating set of rectangles is small enough so that the tail behavior of $\max_{R \in \ell \text{th block}} T(R)$ can be controlled quite precisely by simply using Boole's inequality; see (10) below. This is in contrast to the approximating set of intervals introduced by Rufibach and Walther (2009) for certain multiple testing problems on the line. While that set leads to computationally efficient algorithms, its cardinality is still so large that the control of $\max_{I \in \ell \text{th block}} T(I)$ requires in addition the difficult stochastic control of the increments of $T(I)$ as a process in $I$. In light of the above results, one may surmise that for these and related problems, there typically exists an approximating set that not only allows for statistical optimality and computationally efficient algorithms, but which also obviates the need for the stochastic control of the increments of $T$ when used in conjunction with the block procedure. In particular, it may be possible to recover the optimality results for the inference problems treated in Rufibach and Walther (2009) with the univariate version of the algorithm introduced here.

**4. Algorithm.** It is helpful to use the following notation in this section. The coordinates of the $N$ locations are $(X_1, Y_1), \ldots, (X_N, Y_N)$, and we write $X_{(r)} := X_{(\text{round}(r) \wedge N)}$ for real $r$, where $X_{(1)} \leq \cdots \leq X_{(N)}$ are the order statistics of $X_1, \ldots, X_N$. Each location has a label that is either 0 or 1. Here, is the pseudo-code to enumerate the set of approximating rectangles and to compute the corresponding local test statistics:

Sort the locations $(X_1, Y_1), \ldots, (X_N, Y_N)$ according to the **X**-value.
for $\ell = 3, \ldots, \lfloor \log_2 \frac{N}{2 \log N} \rfloor$ do:
    Set $s := 2^{-\ell}$, $\varepsilon := \ell^{-1/2}/6$.
    for $i = 0, \ldots, \ell$ do:
        for $j = 0, \ldots, \lfloor (\varepsilon s 2^i)^{-1} \rfloor$ do:
            for $k = j + 1, \ldots, j + \lfloor \frac{1}{\varepsilon} \rfloor$ do:
                Extract all locations $(X_p, Y_p)$ for which $X_p$ falls in the interval
                $\mathcal{X}_{jk} := [X_{(j\varepsilon s 2^i N + 1)}, X_{(k\varepsilon s 2^i N)}]$ and denote by
                $N_{jk}$ the number of these locations.
                Sort these extracted $Y_p$ and compute the vector of cumulative
                sums of the labels of the $(X_p, Y_p)$ corresponding to the sorted $Y_p$.
                for $m = 0, \ldots, \lfloor 2^i/\varepsilon \rfloor$ do:
                    for $n = m + 1, \ldots, m + \lfloor 2/\varepsilon \rfloor$ do:
                          Compute the test statistic on the rectangle
                          $\mathcal{X}_{jk} \times [Y_{(m\varepsilon 2^{-i} N_{jk} + 1)}, Y_{(n\varepsilon 2^{-i} N_{jk})}]$, where the order
                          statistics $Y_{(\cdot)}$ are with respect to the extracted $Y_p$.

The running time of the algorithm is almost linear in $N$.

PROPOSITION 1. *The above algorithm runs in $O(N \log^4 N)$ time.*

We illustrate the methodology with an example where 1000 locations are drawn from a mixture of four bivariate normals. The labels are Bernoulli with $p = 0.4$, except in the strip $[x \geq 5]$, where $p = 0.6$, and in the box $[1, 2] \times [3, 5]$, where $p = 0.75$. Critical values for the conventional scan statistic (1) and the blocked scan statistic (3) were computed with 50000 random permutations of the labels. Figure 2 (left) shows all minimal (with respect to inclusion) boxes that are significant at the 5% level using the calibration for the traditional scan statistic. Thus, we are 95% confident that each of the depicted boxes contains a so-called overdensity, that is, somewhere inside the box the probability of success $p$ is larger than outside the box. Figure 2 (right) shows the resulting significant boxes when the calibration for the blocked scan statistic is used. In addition to detecting the small box at $[1,2] \times [3,5]$, the blocked scan statistic also detects the large box $[x \geq 5]$. It was found that this result was a frequent outcome for realizations of this example.

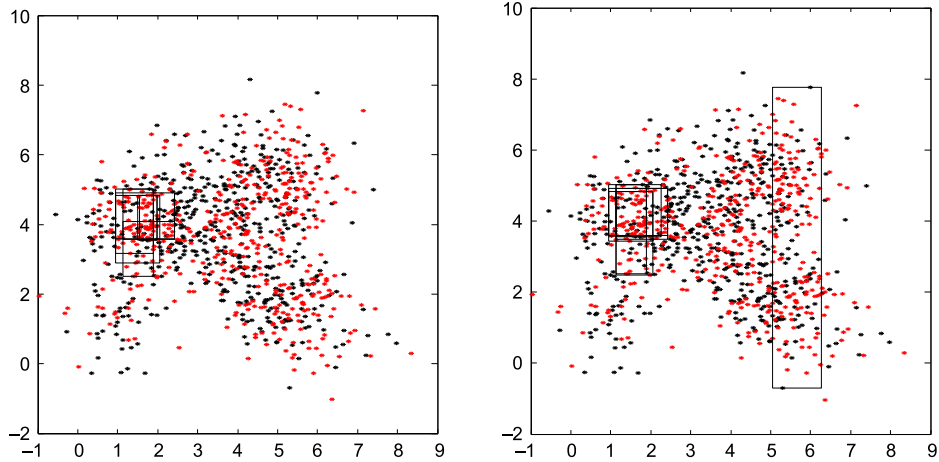

Fig. 2. *Minimal significant rectangles obtained with the conventional calibration for the scan statistic (left) and the blocked scan statistic (right). Locations having label 1 are plotted in red, locations with label 0 are black.*

## 5. Proofs.

PROOF OF THEOREM 1. We parametrize rectangles as follows: $R = (x, x', y, y')$ denotes the rectangle with the vertices $(x, y), (x', y), (x', y')$ and $(x, y')$, where $x, x', y, y' \in [-\infty, \infty]$ and $x < x'$, $y < y'$. We will approximate $\mathcal{R}(s)$ by the finite set $\mathcal{R}_{\mathrm{app}}(s, \varepsilon)$, which for notational simplicity we define for $\varepsilon \in (0, \frac{1}{6})$ by $\mathcal{R}_{\mathrm{app}}(s, 6\varepsilon) := \{R : R = (x_j, x_k, y_m, y_n) := (F_X^{-1}(j\varepsilon s 2^i), F_X^{-1}(k\varepsilon s 2^i), F_Y^{-1}(m\varepsilon 2^{-i}), F_Y^{-1}(n\varepsilon 2^{-i}))$, where $i, j, k, m$ and $n$ are integers with $0 \leq i \leq \lceil \log_2 \frac{1}{s} \rceil$, $0 \leq j \leq \lfloor (\varepsilon s 2^i)^{-1} \rfloor$, $j + 1 \leq k \leq j + \lfloor 1/\varepsilon \rfloor$, $0 \leq m \leq \lfloor 2^i/\varepsilon \rfloor$ and $m + 1 \leq n \leq m + \lfloor 2/\varepsilon \rfloor \}$. Here, $F_X^{-1}$ and $F_Y^{-1}$ denote the quantile functions of the first and second marginals of $F$, respectively, with $F_X^{-1}(p) = -\infty$ for $p < 0$ and $F_X^{-1}(p) = \infty$ for $p > 1$.

Thus,

$$\#\mathcal{R}_{\mathrm{app}}(s, 6\varepsilon) \leq \varepsilon^{-1} 2\varepsilon^{-1} \sum_{i=0}^{\lceil \log_2(1/s) \rceil} ((\varepsilon s)^{-1} 2^{-i} + 1)(\varepsilon^{-1} 2^i + 1)$$

$$\leq 2\varepsilon^{-2} \sum_{i=0}^{\lceil \log_2(1/s) \rceil} \frac{4}{3}(\varepsilon s)^{-1} 2^{-i} \frac{7}{6} \varepsilon^{-1} 2^i$$

$$\leq 4 s^{-1} \varepsilon^{-4} \log_2(4/s).$$

Now, let $R = (x, x'y, y') \in \mathcal{R}(s)$. We will show that there exists a $R' \in \mathcal{R}_{\text{app}}(s, 6\varepsilon)$ with $F(R \triangle R') \leq 6\varepsilon s$ and that one can even arrange that $R' \subset R$. To this end set $i := \lceil \log_2 \frac{F_X([x,x'])}{s} \rceil$. Thus,

(4) $$s 2^{i-1} < F_X([x, x']) \leq s 2^i$$

so the index $i$ is assigned to rectangles whose "length," as measured by $F_X$, lies between $\frac{s}{2} 2^i$ and $s 2^i$. Let $j$ be the smallest integer such that $x_j > x$, let $k$ be the largest integer so that $x_k < x'$, let $m$ be the smallest integer with $y_m > y$, and let $n$ be the largest integer such that $y_n < y'$. It will be shown below that these indices fall in the ranges given in the definition of $\mathcal{R}_{\text{app}}(s, 6\varepsilon)$, hence $R' := (x_j, x_k, y_m, y_n) \in \mathcal{R}_{\text{app}}(s, 6\varepsilon)$, and by definition $R' \subset R$. Further,

(5) $$F(R \triangle R') = F([x, x'] \times ([y, y_m] \cup [y_n, y']))$$
$$+ F(([x, x_j] \cup [x_k, x']) \times [y_m, y_n]).$$

Now, $F_Y([y, y_m]) \leq \varepsilon 2^{-i}$ by the definition of $m$, and the same bound applies to $F_Y([y_n, y'])$. Likewise, both $F_X([x, x_j])$ and $F_X([x_k, x'])$ are not larger than $\varepsilon s 2^i$. Hence, (5) is not larger than

$$2 F_X([x, x']) \varepsilon 2^{-i} + 2 \varepsilon s 2^i F_Y([y, y'])$$
$$= 2\varepsilon(z + s F(R)/z) \qquad \text{where } z := F_X([x, x']) 2^{-i} \in (s/2, s] \text{ by (4)}$$
$$\leq 2\varepsilon(s/2 + 2 F(R)) \qquad \text{since } s/2 < F(R)$$
$$< 6\varepsilon F(R).$$

It remains to show that $i, j, k, m$ and $n$ fall in the ranges given in the definition of $\mathcal{R}_{\text{app}}(s, 6\varepsilon)$: $1 \geq F_X([x, x']) \geq F(R) > s/2$ implies $0 \leq i \leq \lceil \log_2 \frac{1}{s} \rceil$. Clearly, $j \geq 0$. For $\tilde{j} := \lfloor (\varepsilon s 2^i)^{-1} \rfloor$, we have $x_{\tilde{j}} \geq F_X^{-1}(1 - \varepsilon s 2^i) \geq F_X^{-1}(1 - 2\varepsilon F_X([x, x'])) \geq x$ by (4) and as $\varepsilon \leq 1/6$. Hence, $j \leq \tilde{j}$. Next,

$$(k-j)\varepsilon s 2^i = F_X([x_j, x_k]) \begin{cases} \leq F_X([x, x']) \leq s 2^i, \\ \geq F_X([x, x']) - 2\varepsilon s 2^i > 0, \end{cases}$$

by (4) and since $\varepsilon \leq 1/6$. Hence, $1 \leq k - j \leq 1/\varepsilon$. Clearly, $m \geq 0$. Further, $s/2 < F(R) \leq (F_X([x_j, x_k]) + 2\varepsilon s 2^i) \times (F_Y([y_m, \infty)) + \varepsilon 2^{-i}) = (k - j + 2)\varepsilon s 2^i (1 - (m-1)\varepsilon 2^{-i})$. Together with $(k - j + 2)\varepsilon \leq 3$ (see above), this inequality yields $6(m-1)\varepsilon < 6 \times 2^i - 1$, whence $m < 2^i/\varepsilon - 1/(6\varepsilon) + 1 \leq 2^i/\varepsilon$ since $\varepsilon \leq 1/6$. Finally, $s/2 < F(R) \leq (F_X([x_j, x_k]) + 2\varepsilon s 2^i) \times (F_Y([y_m, y_n]) +$ $2\varepsilon 2^{-i}$). With (4) this yields $F_Y([y_m, y_n]) \geq 2^{-i-1}/(1+2\varepsilon) - \varepsilon 2^{-i+1} > 0$ since $\varepsilon \leq 1/6$, hence $n > m$. Likewise, $s \geq F(R) \geq F_X([x_j, x_k]) \times F_Y([y_m, y_n])$ together with (4) gives $s > s2^{i-1}(n-m)\varepsilon 2^{-i}$, hence $n - m < 2/\varepsilon$. □

PROOF OF COROLLARY 1. Define the random collection of rectangles as in (2), where $\mathcal{R}_{\text{app},N}(s, \varepsilon)$ is defined as $\mathcal{R}_{\text{app}}(s, \varepsilon)$ in the proof of Theorem 1 but with $F_X^{-1}$ and $F_Y^{-1}$ replaced by the empirical quantile functions $F_{N,X}^{-1}$ and $F_{N,Y}^{-1}$, respectively. That is, we consider rectangles $R = (x_j, x_k, y_m, y_n)$ such that $(x_j, x_{j+1}]$ has empirical measure $F_{N,X}$ about equal to $\frac{\varepsilon}{6} s 2^i$, likewise for $(x_k, x_{k+1}]$; further $(y_m, y_{m+1}]$ and $(y_n, y_{n+1}]$ have empirical measure $F_{N,Y}$ about equal to $\frac{\varepsilon}{6} 2^{-i}$, where $0 \leq i \leq \lceil \log_2 \frac{1}{s} \rceil$. Then result 2 of Theorem 1 yields

$$\#\mathcal{R}_{\text{app},N} \leq C \sum_{\ell=1}^{\lfloor \log_2(N/(2 \log N)) \rfloor} 2^\ell \ell^2 \log 2^{\ell+1}$$

$$\leq 2C \frac{N}{\log N} \left( \log_2 \frac{N}{2 \log N} \right)^3 \log 2$$

$$\leq 2CN \log_2^2 N.$$

Now, let $R = (x, x', y, y')$ be a rectangle parametrized as in the proof of Theorem 1. Set $\ell := \lfloor \log_2 \frac{1}{F(R)} \rfloor$. Then $3 \leq \ell \leq \lfloor \log_2 \frac{N}{2 \log N} \rfloor$ by the assumptions on $R$. We will construct another deterministic rectangle $\tilde{R} = (\tilde{x}, \tilde{x}', \tilde{y}, \tilde{y}')$ such that for $\gamma = 1/8$:

(A) $F(R \setminus \tilde{R}) \leq (1+\gamma) \frac{F(R)}{\sqrt{\lfloor \log_2(1/F(R)) \rfloor}}$;

(B) $\mathbb{P}(\text{there exists } R' \in \mathcal{R}_{\text{app},N}(2^{-\ell}, \ell^{-1/2}) : \tilde{R} \subset R' \subset R) \geq 1 - 16 \frac{1+\gamma}{\gamma^2 \sqrt{\log N}}$.

The claim of the corollary then follows.

To keep familiar notation set $s := 2^{-\ell}$ and $\varepsilon := \ell^{-1/2}$. Define

$$i := \left\lceil \log_2 \frac{F_X([x, x'])}{s} \right\rceil,$$

so

(6) $$2^{i-1} < \frac{F_X([x, x'])}{s} \leq 2^i.$$

Then $0 \leq i \leq \lceil \log_2 \frac{1}{s} \rceil$ as required in the definition of $\mathcal{R}_{\text{app},N}(s, \varepsilon)$.

To construct $\tilde{R}$, define $\tilde{x}, \tilde{x}', \tilde{y}, \tilde{y}'$ such that $F_X([x, \tilde{x}]) = F_X([\tilde{x}', x']) = (1+\gamma)\varepsilon s 2^i/6$, $F_Y([y, \tilde{y}]) = F_Y([\tilde{y}', y']) = (1+\gamma)\varepsilon 2^{-i}/6$. Then $x < \tilde{x} < \tilde{x}' < x'$ and $y < \tilde{y} < \tilde{y}' < y'$: by (6), $(1+\gamma)\varepsilon s 2^i/6 < F_X([x, x'])/2$. Likewise, the definition of $\ell$ implies $s = 2^{-\ell} < 2F(R)$, hence $(1+\gamma)\varepsilon 2^{-i}/6 \leq \frac{s}{4F_X([x,x'])} = \frac{sF_Y([y,y'])}{4F(R)} < \frac{F_Y([y,y'])}{2}$, which yields $y < \tilde{y} < \tilde{y}' < y'$. Thus, $\tilde{R} \subset R$ and

$$F(R \setminus \tilde{R}) = F([x,x'] \times ([y,\tilde{y}] \cup [\tilde{y}',y']))$$
$$+ F(([x,\tilde{x}] \cup [\tilde{x}',x']) \times [\tilde{y},\tilde{y}'])$$
$$< (1+\gamma)\varepsilon F(R) \qquad \text{as in the proof of Theorem 1,}$$

using $s/2 < F(R)$ by the definition of $\ell$. This establishes (A). Next,

$$\mathbb{P}\left(\text{at most } \frac{\varepsilon}{6}s2^i N \text{ observations in } [x,\tilde{x}]\right)$$
$$= \mathbb{P}(F_{N,X}([x,\tilde{x}]) \leq \varepsilon s 2^i)$$
$$= \mathbb{P}(F_{N,X}([x,\tilde{x}]) - F_X([x,\tilde{x}]) \leq -\gamma\varepsilon s 2^i/6)$$
$$\leq \frac{1+\gamma}{N\gamma^2 \varepsilon s 2^i/6} \qquad \text{by Chebyshev}$$
$$\leq 3\frac{1+\gamma}{\gamma^2 \sqrt{(\log 2)(\log N)}}$$

since $s2^i \geq F_X([x,x']) \geq F(R) \geq 2\frac{\log N}{N}$ by (6) and $\ell \leq \log_2 N$. Hence, with probability at least $1 - 4\frac{1+\gamma}{\gamma^2 \sqrt{\log N}}$ an endpoint $x_j$ of one the rectangles in $\mathcal{R}_{\text{app},N}(2^{-\ell}, \ell^{-1/2})$ falls into $[x,\tilde{x}]$. Analogously, one can show that some $x_k$, $y_m$ and $y_n$ fall into $[\tilde{x}',x']$, $[y,\tilde{y}]$ and $[\tilde{y}',y']$, respectively. Hence, there exists a rectangle $R' \in \mathcal{R}_{\text{app},N}(2^{-\ell}, \ell^{-1/2})$ that satisfies $\tilde{R} \subset R' \subset R$. (B) follows.
□

PROOF OF THEOREM 2(a). To simplify notation we write $(F, R, p, q)$ for $(F^N, R_N, p_N, q_N)$. It will become clear that these parameters may vary with $N$ due to the uniformity of the following results. As usual, $F_N$ will denote the empirical measure pertaining to $F$.

Set $b_N := \varepsilon_N \sqrt{\log \frac{1}{F(R)}} \to \infty$. (a) follows after showing:

(A) $\quad \mathbb{P}_N\Bigg(\text{there exists a rectangle } R' \in \mathcal{R}_{\text{app},N} \text{ with } R' \subset R \text{ and}$

$$T(R') > \log \frac{1}{F(R)} + \sqrt{b_N \log \frac{1}{F(R)}}\Bigg) \to 1.$$

(B) $R'$ belongs to a block $\ell$ whose critical value satisfies

$$q_\ell\left(\frac{\tilde{\alpha}}{\ell^2}\right) \leq \log \frac{1}{F(R)} + 6\log\log \frac{1}{F(R)} + \gamma$$

with $\gamma$ depending on $\alpha$ only.

For the proof of (A), one verifies that the condition $D(F(R),p,q) \geq (2+\varepsilon_N)\log\frac{1}{F(R)}/N$ together with the inequality $\frac{(p-q)^2}{p(1-q)} \leq 1$ for $q < p$ implies $F(R) \geq 2\frac{\log N}{N}$. So the block index $\ell := \lfloor \log_2 \frac{1}{F(R)} \rfloor$ of $R$ satisfies $3 \leq \ell \leq \lfloor \log_2 \frac{N}{2\log N} \rfloor$. Thus, by Corollary 1, with probability converging to 1 there exists $R' \in \mathcal{R}_{\text{app},N}$ such that $R' \subset R$ and

$$F(R') \geq F(R)\left(1 - \frac{9}{8\sqrt{\lfloor \log_2(1/F(R)) \rfloor}}\right)$$
$$\geq F(R)(1-\lambda_N),$$

where $\lambda_N := \frac{2}{3}\sqrt{\frac{b_N}{b_N + \lfloor \log_2(1/F(R)) \rfloor}}$, and the last inequality follows from $3(1+\lfloor \log_2 \frac{1}{F(R)} \rfloor/b_N) \leq (1+3/b_N)\lfloor \log_2 \frac{1}{F(R)} \rfloor \leq \frac{2^8}{3^5}\lfloor \log_2 \frac{1}{F(R)} \rfloor$ for $N$ large enough, using $\lfloor \log_2 \frac{1}{F(R)} \rfloor \geq 3$.

Denote by $\hat{p}$ and $\hat{q}$ the proportion of 1's in $R'$ and $R'^c$, respectively. On the event $\mathcal{A}_N := \{\frac{F_N(R')(1-F_N(R'))}{F(R)(1-F(R))} \geq 1 - \frac{25}{24}\lambda_N, \frac{\hat{p}-\hat{q}}{p-q} \geq 1 - \frac{\lambda_N}{6}, \frac{p}{\hat{p}} \geq 1 - \frac{\lambda_N}{24}, \frac{1-q}{1-\hat{q}} \geq 1 - \frac{\lambda_N}{24}\}$ we have $\hat{q} < \hat{p}$ as $q < p$. The function $l(\hat{p}) := \hat{p}\log\frac{\hat{p}}{\overline{p}} + (1-\hat{p})\log\frac{1-\hat{p}}{1-\overline{p}}$ satisfies $l(\overline{p}) = 0$, $l'(\overline{p}) = 0$, and $l''(\xi) = \xi^{-1}(1-\xi)^{-1} \geq \hat{p}^{-1}(1-\hat{q})^{-1}$ for $\xi \in [\hat{q},\hat{p}]$. Thus, Taylor's theorem gives on $\mathcal{A}_N$

$$\begin{aligned}
T(R') &= (\#R')l(\hat{p}) + (N - \#R')l(\hat{q}) \\
&\geq (\#R')\frac{(\hat{p}-\overline{p})^2}{2\hat{p}(1-\hat{q})} + (N - \#R')\frac{(\hat{q}-\overline{p})^2}{2\hat{p}(1-\hat{q})} \\
&= \frac{(\#R')(N-\#R')}{N}\frac{(\hat{p}-\hat{q})^2}{2\hat{p}(1-\hat{q})} \\
&\geq NF(R)(1-F(R))\frac{(p-q)^2}{2p(1-q)}\left(1 - \frac{25}{24}\lambda_N\right)\left(1 - \frac{\lambda_N}{6}\right)^2\left(1 - \frac{\lambda_N}{24}\right)^2 \\
&\geq \left(1 - \frac{35}{24}\lambda_N\right)ND(F(R),p,q)/2 \\
&\geq \log\frac{1}{F(R)} + \frac{1}{2}\Biggl(\left(1 - \frac{35}{24}\lambda_N\right)b_N - \frac{35}{12}\lambda_N\sqrt{\log\frac{1}{F(R)}}\Biggr)\sqrt{\log\frac{1}{F(R)}} \\
&\geq \log\frac{1}{F(R)} + \frac{1}{72}(b_N - 140\sqrt{b_N})\sqrt{\log\frac{1}{F(R)}} \qquad \text{since } \lambda_N \leq 2/3.
\end{aligned}$$

(A) follows once we show that $\mathbb{P}_N(\mathcal{A}_N) \to 1$. As for the first event in $\mathcal{A}_N$, the proof of Corollary 1 provided a deterministic rectangle $\tilde{R}$ with $F(\tilde{R}) \geq (1-\lambda_N)F(R) \geq \frac{2}{3}\log N/N$ and $\mathbb{P}_N(\tilde{R} \subset R' \subset R) \to 1$. Chebyshev's inequality gives for $c = 1/24$:

$$\mathbb{P}_N\left(\left|\frac{F_N(\tilde{R})}{F(\tilde{R})} - 1\right| > c\lambda_N\right)$$
$$\leq \frac{F(\tilde{R})}{N(F(\tilde{R}))^2 c^2 \lambda_N^2} \leq \frac{3}{2(\log N)\lambda_N^2 c^2}$$
$$\leq \frac{27}{4c^2 \min(b_N \log 2, \log N)} \quad \text{by (8)}$$

and the same bound holds for $|\frac{F_N(R)}{F(R)} - 1|$. But $\tilde{R} \subset R' \subset R$, $F_N(\tilde{R})/F(\tilde{R}) \geq 1 - \lambda_N/24$ and $F_N(R)/F(R) \leq 1 + \lambda_N/24$ imply

$$(7) \quad \frac{F_N(R')}{F(R)} \begin{cases} \geq \frac{F_N(\tilde{R})}{F(\tilde{R})}\frac{F(\tilde{R})}{F(R)} \geq \left(1 - \frac{\lambda_N}{24}\right)(1-\lambda_N) \geq 1 - \frac{25}{24}\lambda_N, \\ \leq \frac{F_N(R)}{F(R)} \leq 1 + \frac{\lambda_N}{24}. \end{cases}$$

This entails the first event in $\mathcal{A}_N$ due to the inequality $\frac{y(1-y)}{x(1-x)} \geq \min(\frac{y}{x}, 2 - \frac{y}{x})$ for $x, y \in (0, 1/2)$.

For the other events in $\mathcal{A}_N$, note that given the locations $\mathbf{X} = (X_1, \ldots, X_N)$, $\hat{p}$ and $\hat{q}$ are independent with $\hat{p} \sim \text{bin}(\#R', p)/\#R'$, while $\hat{q}$ is an average of $N - \#R'$ independent Bernoulli random variables that have probability of success equal to $p$ for the $\#R - \#R'$ locations in $R \setminus R'$ and $q$ for the $N - \#R$ locations in $R^c$. Hence,

$$\mathbb{E}_N(\hat{q}|\mathbf{X}) = \frac{(\#R - \#R')p + (N - \#R)q}{N - \#R'}$$
$$= q + \frac{F_N(R) - F_N(R')}{1 - F_N(R')}(p - q)$$
$$\leq q + \frac{13/(12 \cdot 8)\lambda_N}{1 - 37/(8 \cdot 36)}(p - q) \quad \text{on (7) as } \lambda_N \leq \frac{2}{3}, F(R) \leq \frac{1}{8}$$
$$< q + \frac{3}{19}\lambda_N(p - q).$$

Thus, on (7)

$$\mathbb{P}_N\left(\frac{\hat{p} - \hat{q}}{p - q} < 1 - \frac{\lambda_N}{6}\bigg|\mathbf{X}\right)$$
$$\leq \mathbb{P}_N\left(\hat{p} - \hat{q} - \mathbb{E}_N(\hat{p} - \hat{q}) < \left(\frac{3}{19} - \frac{1}{6}\right)\lambda_N(p - q)\bigg|\mathbf{X}\right)$$

$$\leq \frac{p(1-q)/\#R' + p(1-q)/(N-\#R')}{114^{-2}\lambda_N^2(p-q)^2} \quad \text{as } q < p$$

$$= \frac{114^2}{\lambda_N^2 N D(F(R), p, q)} \cdot \frac{F(R)(1-F(R))}{F_N(R')(1-F_N(R'))}$$

$$\leq \frac{114^2}{\lambda_N^2 (2 + b_N (\log(1/F(R)))^{-1/2}) \log(1/F(R))} \cdot \left(1 - \frac{25}{24}\lambda_N\right)^{-1}$$

$$\leq \frac{9 \cdot 114^2 \cdot 2}{b_N} \to 0 \quad \text{by (9) and as } \lambda_N \leq \frac{2}{3}.$$

The other two events in $\mathcal{A}_N$ obtain similarly. Above, we used the following properties of $\lambda_N$:

(8) $\qquad\qquad\qquad\qquad \lambda_N^2 \log N \geq \frac{2}{9} \min(b_N \log 2, \log N),$

(9) $\qquad \lambda_N^2 \left(2 + \frac{b_N}{\sqrt{\log(1/F(R))}}\right) \log \frac{1}{F(R)} \geq 2b_N/9.$

For proof of those inequalities, note that $\lambda_N^2 \log_2 \frac{1}{F(R)} \geq \frac{4}{9} \cdot \frac{b_N \log_2(1/F(R))}{b_N + \log_2(1/F(R))} \geq \frac{2}{9} \min(b_N, \log_2 \frac{1}{F(R)})$. Now (8) follows since $F(R) \geq 2\log N/N$ implies $\log N / \log_2 \frac{1}{F(R)} \geq \log 2$. Applying the above inequality to the LHS of (9), yields the lower bound $\frac{2}{9} \min(b_N 2/\log_2(e), b_N \sqrt{\log \frac{1}{F(R)}}) \geq \frac{2}{9} b_N$.

For the proof of (B), note that by the construction of $\mathcal{R}_{\text{app},N}$ in the proof of Corollary 1, the rectangle $R'$ belongs to $\mathcal{R}_{\text{app},N}(2^{-\ell}, \ell^{-1/2})$ where $\ell := \lfloor \log_2 \frac{1}{F(R)} \rfloor$. Hence, result 2 of Theorem 1 yields $\#\mathcal{R}_{\text{app},N}(2^{-\ell}, \ell^{-1/2}) \leq K \frac{1}{F(R)} (\log_2 \frac{1}{F(R)})^3$ for some universal constant $K$. Thus,

$$\mathbb{P}_0\left(\max_{\tilde{R} \in \mathcal{R}_{\text{app},N}(2^{-\ell}, \ell^{-1/2})} T(\tilde{R}) \geq \log \frac{1}{F(R)} + 6 \log\log \frac{1}{F(R)} + \gamma\right)$$

(10) $\qquad \leq \frac{K}{F(R)} \left(\log_2 \frac{1}{F(R)}\right)^3$

$$\times \max_{\tilde{R} \in \mathcal{R}_{\text{app},N}(2^{-\ell}, \ell^{-1/2})} \mathbb{P}_0\left(T(\tilde{R}) \geq \log \frac{1}{F(R)} + 6 \log\log \frac{1}{F(R)} + \gamma\right),$$

where $\mathbb{P}_0$ denotes the null distribution, i.e., the permutation distribution, conditional on the $N$ locations and on $\overline{p} = \frac{\text{total no. of 1's}}{N}$. Thus, under $\mathbb{P}_0$, the number of 1's in a rectangle $\tilde{R}$ follows the hypergeometric distribution where $\#\tilde{R}$ labels are drawn out of $N$, of which $\overline{p}N$ are 1's. Theorem 4 implies

$$\mathbb{P}_0\left(T(\tilde{R}) \geq \log \frac{1}{F(R)} + 6 \log\log \frac{1}{F(R)} + \gamma\right)$$

$$\leq 2C\left(\log\frac{1}{F(R)} + 6\log\log\frac{1}{F(R)} + \gamma + 2\right)$$
$$\times F(R)\left(\log\frac{1}{F(R)}\right)^{-6}\exp(-\gamma)$$
$$\leq \frac{6\alpha}{K\pi^2}F(R)\left(\log_2\frac{1}{F(R)}\right)^{-5}$$

for $\gamma$ large enough, depending on $\alpha$ only. Thus, (10) is not larger than $\frac{6\alpha}{\pi^2} \times (\log_2 \frac{1}{F(R)})^{-2} \leq \frac{6\alpha}{\pi^2 \ell^2}$ by the definition of $\ell$. Now, (B) follows once it is shown that

(11) $$\mathbb{P}_0\left(\max_{\tilde{R}\in\mathcal{R}_{\mathrm{app},N}(2^{-\ell},\ell^{-1/2})}T(\tilde{R}) \geq q_\ell\left(\frac{\tilde\alpha}{\ell^2}\right)\right) \geq \frac{6\alpha}{\pi^2\ell^2}.$$

But by the definition of $q_\ell(\cdot)$, the probability in (3) is not larger than $\sum_{\ell\geq 1}\tilde\alpha/\ell^2 \leq \tilde\alpha\pi^2/6$, hence $\tilde\alpha \geq 6\alpha/\pi^2$ by the definition of $\tilde\alpha$. Now (11) follows from the definition of $q_\ell(\cdot)$. □

PROOF OF THEOREM 2(b1, b2). The idea of the proof of (b1) is classical, see, e.g., Lepski and Tsybakov (2000). Given the prescribed sequence of values $\{F^N(R_N)\}$ and $\{\varepsilon_N\}$, partition $\mathbf{R}^2$ into rectangles $\tilde R_1^N, \tilde R_2^N, \ldots$, such that $F^N(\tilde R_j^N) = F^N(R_N)$ for $j=1,\ldots,\lfloor\frac{1}{F^N(R_N)}\rfloor$. This is feasible since $F^N$ is continuous, e.g., by partitioning one axis into intervals. Set $q = q_N := 1/2$ and $p = p_N := q + \sqrt{\frac{1}{2NF(R)(1-F(R))}\log\frac{1}{F(R)}}(1-\varepsilon_N/8)$, where for notational simplicity, we write $F(R)$ for $F^N(R_N)$ and also drop the index $N$ from $F^N, p_N, q_N, \tilde R_j^N$ in the following. Without loss of generality, we may assume $\varepsilon_N < 8$. Thus, $q < p$ and for $j=1,\ldots,\lfloor\frac{1}{F(R)}\rfloor$:

$$D(F(\tilde R_j), p, q)$$
$$= \frac{\log(1/F(R))(1-\varepsilon_n/8)^2}{Np}$$
$$\geq \frac{\log(1/F(R))(1-\varepsilon_n/4)/N}{1/2 + (\log(1/F(R)))^{-1/2}} \quad \text{by (12)}$$
$$\geq (2-\varepsilon_N/2)\left(1 - 2\left(\log\frac{1}{F(R)}\right)^{-1/2}\right)\frac{\log(1/F(R))}{N}$$
$$\geq (2-\varepsilon_N)\frac{\log(1/F(R))}{N}$$

for $N$ large enough, as $\varepsilon_N \sqrt{\log \frac{1}{F(R)}} \to \infty$. We used

$$(12) \quad NF(R)(1-F(R)) \geq (1+o(1))\log^2 N \geq (1+o(1))\log^2 \frac{1}{F(R)},$$

which is a consequence of the assumptions on $F(R)$ stated in the theorem.

Denote by $X_i$ the location and by $Y_i$ the Bernoulli random variable that gives the label of the $i$th observation, $i=1,\ldots,N$. Denote by $\mathbb{P}_{N,0}$ the model where the $X_i$ are i.i.d. $F$ and all the $Y_i$ have probability of success $q$, while we define $\mathbb{P}_{N,j}$ to be the model where instead $Y_i$ has probability of success $p$ if $X_i \in \tilde{R}_j$ and $q$ otherwise, $i=1,\ldots,N$. Thus, $\mathbb{P}_{N,0}$ belongs to $H_0$. Define the likelihood ratio $L_{N,j}(\mathbf{X},\mathbf{Y}) := \prod_{i=1}^N f_{N,j}(X_i,Y_i)$, where

$$f_{N,j}(X_i,Y_i) := \begin{cases} \frac{p}{q}1(Y_i=1) + \frac{1-p}{1-q}1(Y_i=0), & \text{if } X_i \in \tilde{R}_j, \\ 1, & \text{otherwise,} \end{cases}$$

$j=1,\ldots,\lfloor \frac{1}{F(R)} \rfloor$. Hence, if $\phi_N(\mathbf{X},\mathbf{Y})$ is any level $\alpha$ test that depends on the locations $\mathbf{X}$ and the labels $\mathbf{Y}$, then by conditioning on $\mathbf{X}$ one verifies $\mathbb{E}_{N,j}\phi_N(\mathbf{X},\mathbf{Y}) = \mathbb{E}_{N,0}\phi_N(\mathbf{X},\mathbf{Y})L_{N,j}(\mathbf{X},\mathbf{Y})$. We will show that

$$(13) \quad \mathbb{E}_{N,0}\left| \left\lfloor \frac{1}{F(R)} \right\rfloor^{-1} \sum_{j=1}^{\lfloor 1/F(R) \rfloor} L_{N,j}(\mathbf{X},\mathbf{Y}) - 1 \right| \to 0.$$

Then

$$\min_{j=1,\ldots,\lfloor 1/F(R) \rfloor} \mathbb{E}_{N,j}\phi_N(\mathbf{X},\mathbf{Y}) - \alpha$$

$$\leq \left\lfloor \frac{1}{F(R)} \right\rfloor^{-1} \sum_{j=1}^{\lfloor 1/F(R) \rfloor} \mathbb{E}_{N,j}\phi_N(\mathbf{X},\mathbf{Y}) - \alpha$$

$$= \mathbb{E}_{N,0}\left( \left\lfloor \frac{1}{F(R)} \right\rfloor^{-1} \sum_{j=1}^{\lfloor 1/F(R) \rfloor} L_{N,j}(\mathbf{X},\mathbf{Y}) - 1 \right)\phi_N(\mathbf{X},\mathbf{Y}) + o(1)$$

$$\leq \mathbb{E}_{N,0}\left| \left\lfloor \frac{1}{F(R)} \right\rfloor^{-1} \sum_{j=1}^{\lfloor 1/F(R) \rfloor} L_{N,j}(\mathbf{X},\mathbf{Y}) - 1 \right| + o(1)$$

$$= o(1)$$

and the claim of (b1) follows. [Note that one can even allow $\phi_N(\mathbf{X},\mathbf{Y})$ to depend on $F^N(R_N)$. Further, the continuity assumption on $F^N$ was only used to allow for the above partition of $\mathbf{R}^2$ into rectangles and can be relaxed accordingly.] To prove (13), note that conditional on $\mathbf{X}$ the $L_{N,j}(\mathbf{X},\mathbf{Y})$ are independent since $L_{N,j}(\mathbf{X}, \mathbf{Y})$ is a function of only those $Y_i$ for which $X_i \in \tilde{R}_j$. Further, one verifies $\mathbb{E}_{N,0} L_{N,j}(\mathbf{X}, \mathbf{Y}) = \mathbb{E}_{N,0}(L_{N,j}(\mathbf{X}, \mathbf{Y})|\mathbf{X}) = 1$. Thus, we can proceed similarly as in the proof of Lemma 7.4 in Dümbgen and Walther (2008) and obtain (13) once we show that

$$\max\left(\left|\frac{p}{q} - 1\right|, \left|\frac{1-p}{1-q} - 1\right|\right) \leq C\left(\log\left\lfloor \frac{1}{F(R)} \right\rfloor\right)^{-1/2}$$

(14)

for some constant $C$,

(15) $\quad \sqrt{\log\left\lfloor \dfrac{1}{F(R)} \right\rfloor} \left(1 - \dfrac{N \mathbb{E}_{N,0}(f_{N,1} - 1)^2}{2\log\lfloor 1/F(R) \rfloor}\right) \to \infty.$

Now $|\frac{p}{q} - 1| = 2\sqrt{\frac{1}{2NF(R)(1-F(R))} \log \frac{1}{F(R)}}(1 - \varepsilon_N/8) \leq 2(\log\lfloor \frac{1}{F(R)} \rfloor)^{-1/2}$ by (12) and the same bound obtains for $|\frac{1-p}{1-q} - 1|$, proving (14). Finally,

$$\mathbb{E}_{N,0}(f_{N,1}(X_1, Y_1) - 1)^2$$
$$= \mathbb{E}_{N,0}\left(\left(\frac{p}{q} 1(Y_i = 1) + \frac{1-p}{1-q} 1(Y_i = 0) - 1\right)^2 \bigg| X_1 \in \tilde{R}_1\right)$$
$$\times \mathbb{P}_{N,0}(X_1 \in \tilde{R}_1)$$
$$= \frac{2\log(1/F(R))}{NF(R)(1-F(R))}(1 - \varepsilon_N/8)^2 F(R).$$

Together with $\frac{\log x}{\log\lfloor x \rfloor} \leq 1 + \frac{1}{\log\lfloor x \rfloor}$ for $x \geq 2$ one sees that the expression in (15) is not smaller than

$$\sqrt{\log\left\lfloor \frac{1}{F(R)} \right\rfloor}\left[1 - \frac{(1 + 1/\log\lfloor 1/F(R) \rfloor)(1 - \varepsilon_N/8)}{1 - F(R)}\right]$$
$$\geq \sqrt{\log\left\lfloor \frac{1}{F(R)} \right\rfloor}\left[\varepsilon_N/8 - \frac{1}{\log\lfloor 1/F(R) \rfloor} - F(R)\right]$$
$$\bigg/ (1 - F(R))$$
$$\to \infty$$

as $F(R) \to 0$ and $\varepsilon_N \sqrt{\log \frac{1}{F(R)}} \to \infty$, completing the proof of (b1). The bounds on $F^N(R_N)$ and $b_N$ given in the statement of (b2) guarantee that there exists $p_N$ and $q_N$ such that $D(F^N(R_N), p_N, q_N) \geq b_N/N$, e.g., take $p_N = 1$, $q_N = 0$. Then the claim obtains with a contiguity argument similar as in the proof of Theorem 4.1(c) in Dümbgen and Walther (2008). $\square$

PROOF OF THEOREM 3. Part (a) continues to hold as intervals on the axes are special cases of axis-parallel rectangles. Parts (b1) and (b2) continue to hold as their proofs do not depend on the dimensionality of the space. In fact, the proof of (b1) already uses a univariate partitioning of one axis into $\lfloor \frac{1}{F^N(R_N)} \rfloor$ intervals, and the rest of the proof of (b1) goes through verbatim. □

PROOF OF THEOREM 4. Let $k, x \geq 0$ be integers with $x + k \leq \min(n, R)$. Then

$$\frac{\mathbb{P}(X = x+k)}{\mathbb{P}(X = x)} = \prod_{i=1}^{k} \frac{(R-x-k+i)(n-x-k+i)}{(x+i)(N-R-n+x+i)}$$

is nonincreasing in $x$. Hence, for $x \geq \lceil m \rceil$

$$\mathbb{P}(X \geq x) = \frac{\mathbb{P}(X=x)}{\mathbb{P}(X=\lceil m \rceil)} \sum_{k \geq 0} \mathbb{P}(X = \lceil m \rceil) \frac{\mathbb{P}(X=x+k)}{\mathbb{P}(X=x)}$$

$$\leq \frac{\mathbb{P}(X=x)}{\mathbb{P}(X=\lceil m \rceil)} \sum_{k \geq 0} \mathbb{P}(X = \lceil m \rceil) \frac{\mathbb{P}(X=\lceil m \rceil + k)}{\mathbb{P}(X=\lceil m \rceil)}$$

$$\leq \frac{\mathbb{P}(X=x)}{\mathbb{P}(X=\lceil m \rceil)}.$$

The connection between this hypergeometric probability and the log likelihood ratio statistic $L$ obtains by applying Stirling's formula and collecting terms: the upper and lower bounds for Stirling's formula in Feller [(1968), page 54] yield

$$\log \frac{\mathbb{P}(X=x)}{\mathbb{P}(X=\lceil m \rceil)}$$
$$\leq -L(x) + L(\lceil m \rceil)$$
$$+ \frac{1}{2} \log \frac{\lceil m \rceil (R - \lceil m \rceil)(n - \lceil m \rceil)(N - R - n + \lceil m \rceil)}{x(R-x)(n-x)(N-R-n+x)}$$
$$+ \frac{1}{12\overline{p}(1-\overline{p})} \left( \frac{1}{n} + \frac{1}{N-n} \right).$$

Set $L(n, \hat{p}) := n(\hat{p} \log \frac{\hat{p}}{\overline{p}} + (1-\hat{p}) \log \frac{1-\hat{p}}{1-\overline{p}})$. Using $\log \frac{b}{a} \leq \frac{b-a}{a}$ for $0 < a < b$ and Taylor's theorem, respectively, one finds

$$n \frac{(\hat{p} - \overline{p})^2}{\overline{p}(1-\overline{p})} \geq L(n, \hat{p}) \geq \begin{cases} n \dfrac{(\hat{p} - \overline{p})^2}{2(1-\overline{p})}, & \text{if } \hat{p} \geq \overline{p}, \\ n \dfrac{(\hat{p} - \overline{p})^2}{2\overline{p}}, & \text{if } \hat{p} \leq \overline{p}, \end{cases}$$

which implies $L(\lceil m \rceil) \leq \frac{1}{\overline{p}(1-\overline{p})}(\frac{1}{n} + \frac{1}{N-n})$ and for $\hat{p} \in (\overline{p}, \frac{n-1}{n}]$:

$$L(n,\hat{p}) + 1 \geq \frac{1-\overline{p}}{4(1-\hat{p})} \tag{16}$$

[distinguish the cases $\hat{p} \lesseqgtr (3+\overline{p})/4$], as well as for $\hat{q} \in [\frac{1}{N-n}, \overline{p})$:

$$L(N-n,\hat{q}) + 1 \geq \frac{\overline{p}}{4\hat{q}}. \tag{17}$$

First, consider the case $\lceil m \rceil \leq x < \min(n, R)$. Then (17) gives

$$\frac{\lceil m \rceil (R - \lceil m \rceil)}{x(R-x)} \leq \frac{R-m}{R-x} = \frac{\overline{p}}{\hat{q}} \leq 4L(N-n,\hat{q}) + 4$$

and analogously (16) implies

$$\frac{(n - \lceil m \rceil)(N - R - n + \lceil m \rceil)}{(n-x)(N-R-n+x)} \leq 4L(n,\hat{p}) + 4.$$

The first inequality of the theorem now follows from the arithmetic–geometric means inequality.

The case $x = \min(n, R)$ is treated similarly. For example, if $x = n < R$ then $\log \mathbb{P}(X \geq x) = \log \mathbb{P}(X = x) \leq -L(X) + \frac{1}{2}\log(\frac{R(N-n)}{(R-n)N}) + \frac{1}{12}(\frac{1}{n} + \frac{1}{N-n})$ and (17) gives $\frac{R(N-n)}{(R-n)N} = \frac{\overline{p}}{\hat{q}} \leq 4L(N-n,\hat{q}) + 4$, which yields the claimed inequality.

The second inequality of the theorem obtains analogously. The third inequality follows from the first two because the function $x \to L(x)$ is strictly decreasing for $x < m$ and strictly increasing for $x > m$. $\square$

PROOF OF PROPOSITION 1. Sorting the data according to the $X$-coordinate requires $O(N \log N)$ steps. Note that the test statistic inside the $n$-loop can be computed in constant time: the rectangle $\mathcal{X}_{jk} \times [Y_{(a)}, Y_{(b)}]$ contains $\text{round}(b) - \text{round}(a) + 1$ locations, and the number of their labels that equal 1 is just the cumulative sum vector of the labels evaluated at index $\text{round}(b)$ minus the vector evaluated at $\text{round}(a-1)$. These two quantities are sufficient to calculate the test statistic once the overall number of locations $N$ and the overall number of 1's is known. Thus, there are $O(1/\varepsilon)$ steps for the $n$-loop, and hence $O(2^i/\varepsilon^2)$ for the $m$-loop. Inside the $k$-loop the number of steps required to extract the $N_{jk}$ locations $(X_p, Y_p)$, to sort the corresponding $Y$-values, and to compute the cumulative sum is dominated by the sorting, which requires $O(N_{jk} \log N_{jk})$ steps. (Note that presorting the locations according to their $X$-coordinate allows an efficient extraction.)

Thus, the total number of steps in the algorithm is bounded by

$$O(N \log N) + \sum_{\ell=3}^{\log_2(N/(2\log N))} \sum_{i=0}^{\ell} \sum_{j=0}^{\ell^{1/2} 2^\ell 2^{-i}} \sum_{k=j+1}^{j+\ell^{1/2}} (O(N_{jk} \log N_{jk}) + O(2^i \ell)).$$

By definition, $N_{jk} \leq (k-j)\varepsilon s 2^i N \leq 2^{i-\ell} N$. Thus, the above sum is not larger than

$$O(N \log N) + C \sum_{\ell=3}^{\log_2(N/(2\log N))} \sum_{i=0}^{\ell} \ell 2^\ell 2^{-i} (2^{i-\ell} N \log N + 2^i \ell)$$

$$\leq CN \log N \sum_{\ell=3}^{\log_2(N/(2\log N))} \ell^2$$

$$\leq CN(\log N)^4,$$

where the constant $C$ may change from line to line. $\square$

**Acknowledgments.** I would like to thank one referee for the Derbeko, El-Yaniv and Meir (2004) reference and another referee for the Rohde (2009) reference.

Department of Statistics
Stanford University
390 Serra Mall
Stanford, California 94305
USA
E-mail: gwalther@stanford.edu